\documentclass[12pt]{article}
\usepackage[utf8]{inputenc}
\usepackage{amsmath} 
\usepackage{graphicx} 
\usepackage{makeidx}	 
\usepackage{latexsym}
\usepackage{amssymb}
\usepackage{amsfonts}
\usepackage{amsthm}
\newfont{\rams}{msbm10 scaled\magstep1} 

\newtheorem{teo}{Theorem}[section]
\newtheorem{cor}[teo]{Corollary}
\newtheorem{lem}[teo]{Lemma}
\newtheorem{pro}[teo]{Proposition}

\theoremstyle{remark}
\newtheorem{rem}[teo]{Remark}
\newtheorem{oss}[teo]{Observation}

\newtheorem*{que}{Question}

\theoremstyle{definition} 
\newtheorem{defn}[teo]{Definition}

\newenvironment{dimo}%
	{\noindent \underline{\emph{Proof}}  }%
	{\hfill $\Box$ \newline }

	{\underline{\emph{Example}} \newline }%
	{\hfill $\Box$ \newline \newline}
	
\title{Green's Hyperplane Restriction Theorem:\\ an extension to modules}
\author{Ornella Greco}
\date{\empty}
\begin{document}

\maketitle
\begin{abstract}
\noindent
In this paper, we prove a generalization of Green's Hyperplane Restriction Theorem to the case of modules over the polynomial
ring, providing in particular  an upper bound for the Hilbert function of the general linear restriction of a module $M$ in a degree $d$
by the corresponding Hilbert function of a lexicographic module.
\end{abstract}

\section{Introduction}
The extremal properties of the Hilbert function constitute a relevant topic in commutative algebra and, within this topic, an important role  is 
played by lexicographic ideals and binomial representation of integers.\\
Many results about extremal behavior of Hilbert functions have been proved. In particular, Macaulay's Theorem (see \cite{BrunsHerzog, robb, macaulay}),
which characterizes  the possible Hilbert functions of homogeneous $k$-algebras;
or, the analogous result for exterior algebra, proved by Kruskal and Katona (see \cite{katona,kru}), which also provides a 
characterization of $f$-vectors of simplicial complexes.\\
More recently, in \cite{hulett2}, Hulett proved an extension of Macaulay's Theorem to modules over the symmetric algebra.\\
Another classical, but slightly recent, result in this topic is the Hyperplane Restriction Theorem (HRT), given by  Green in 
\cite{green},
which gives a bound for the codimension of the generic linear restriction of a vector space generated 
in a certain degree by the  codimension of a $lex$-segment space with same degree and dimension. This result was also used by Green to 
give another,
less technical, proof of Macaulay's Theorem.  \\
It has been also been applied by several authors 
(see for instance \cite{ahn,boij,migliore}) to get some results about level and Gorenstein algebras, 
with focus on the weak Lefschetz property.\\

In this paper, we will focus on  Green's theorem, trying to follow Hulett's path in generalizing these results to
modules over the polynomial ring.\\
One generalization in this direction has already been done by Gasharov (\cite{gasha}), but we will give an extension of Green's 
theorem
that provides a bound which is actually achieved by lexicographic modules.\\

We will first give the tools to state Green's Hyperplane Restriction Theorem: namely we introduce the binomial 
expansion of a positive integer, also called Macaulay representation, then, the concept of generic linear forms.
Afterwards, we will state  Green's theorem.\\
In the next section, we will extend some of these definitions to the submodule of a finitely generated free module, in particular we will define a 
monomial order induced by the $deglex$ ordering, and will provide the concept of lexicographic module, which is 
a particular class of monomial modules.\\
We first try to extend  Green's theorem in the simplest case of a rank
$2$ module, by proving an inequality of Macaulay representations. Later we prove by induction a new inequality which
imply the main theorem, namely the extension of the Hyperplane Restriction 
Theorem to the case of modules over the polynomial ring.

As a consequence of our theorem, we will derive a result which constitute another version of the theorem "after scaling",
which means that we divide by the respective Hilbert function of the polynomial ring.

Finally, we apply our main theorem to level algebras, by giving some conditions in which the bound given by the extension to modules of the Green's HRT is given than the one obtained trough the Green's HRT.

\section{Green's Theorem}
Let $S=k[x_1, \dots ,x_n]$, where $k$ is an infinite field, and let $S$ be standard graded. Let us fix the $deglex$  monomial 
ordering on $S$
with $x_1>x_2>\cdots >x_n$, namely
\[
  \mathbf{x}^\mathbf{a}>_{deglex}  \mathbf{x}^\mathbf{b} \ \Leftrightarrow \ | \mathbf{a} |  > |  \mathbf{b} | \
  \textrm{or} \ \ |  \mathbf{a} | =|  \mathbf{b} |  \ \textrm{and}\  \mathbf{a}>_{lex} \mathbf{b},
\] (see \cite{robb2}). Notice that in a homogeneous component $S_d$ the monomial orderings, $deglex$ and $lex$, coincide.

\begin{defn} Let $a, d \in \mathbb{N}$, then the \emph{$d$-th Macaulay representation of $a$} (also called the binomial 
representation of $a$ in base $d$) is the unique way to write
$$
a=\binom{a_d}{d}+\binom{a_{d-1}}{d-1}+\cdots + \binom{a_\delta}{\delta},
$$
where $a_d>a_{d-1}>\cdots >a_\delta \geq \delta$
and $\delta=\mathrm{min}\{
i \ | \ a_i \geq i\}$.\\
Let us set $\binom{c}{d}=0$ whenever $c<d$, then, given the $d$-th Macaulay representation of $a$, 
let us define the integer 
$$a_{\langle d\rangle}=\binom{a_d-1}{d}+\binom{a_{d-1}-1}{d-1}+\cdots + \binom{a_\delta-1}{\delta}.
$$
\end{defn}

\noindent
Sometimes it will be convenient to use the extended $d$-th Macaulay representation, i.e. 
$\binom{a_d}{d}+\binom{a_{d-1}}{d-1}+\cdots + \binom{a_1}{1}$, where $a_d>a_{d-1}>\cdots>a_1 \geq 0$. \\
Moreover, 
since $\binom{a}{b}=0$ if $a<b$, we will set the numerators of such binomial coefficients in the  Macaulay representation equal to zero.

\begin{oss}
	If $a=\binom{a_d}{d}+\binom{a_{d-1}}{d-1}+\cdots + \binom{a_1}{1}$ and $b=\binom{b_d}{d}+
	\binom{b_{d-1}}{d-1}+\cdots + \binom{b_1}{1}$,
	 $a\geq b$ if  and only if $(a_d, \dots, a_1)\geq_{lex} (b_d, \dots , b_1)$.
\end{oss}
\begin{defn}
Let $d \in \mathbb{N}$, a \emph{$lex$-segment} is a set constituted by all monomials in $S_d$ 
lexicographically larger or equal than $f$, for some $f\in \mathrm{Mon}(S_d)$. A $k$-vector subspace $V$ of $S_d$
is a \emph{$lex$-segment space} if $\mathrm{Mon}(S_d)\cap V$ is a $k$-basis of $V$ and a $lex$-segment.
\end{defn}

\begin{defn}
A graded monomial ideal $I$ in $S$ is said to be \emph{lexicographic} if, for all degrees $d$, the homogeneous component
$I_d$ is a $lex$-segment space.
\end{defn}
\begin{defn}
	We say that a property $\mathcal{P}$ holds for a generic linear form $\ell$ if there is a non-empty Zariski
	open set $\mathcal{U}\subseteq S_1$ such
	that $\mathcal{P}$ holds for all $\ell \in \mathcal{U}$.
\end{defn}
If $V\subseteq S_d$ a $lex$-segment space, we denote by $V_\ell$ the image of $V$ in $S/(\ell)$.
\begin{pro}\cite[Proposition 5.5.23]{robb}
	Let $k$ be an infinite field, $d\in \mathbb{N}$, $V\subseteq S_d$ a $lex$-segment space. Then:
	\begin{enumerate}
		\item For a generic linear form $\ell \in S_1$, there is a homogeneous linear change of coordinates $\phi:S \rightarrow S$
		such that $\phi(V)=V$ and $\phi(\ell)=x_n$;
		\item For a generic linear form $\ell \in S_1$, we have
				\[
					\mathrm{codim}_k(V_\ell)=\mathrm{codim}_k(V_{x_n})=\mathrm{codim}_k(V)_{\langle d\rangle}.
				\]
	\end{enumerate}
	
\end{pro}
\noindent
So, after a general change of coordinates, we may assume that the generic linear form is $x_n$.\\

\noindent
From now on, we will denote by $H(M,d)$ the Hilbert function of the $S$-module $M$  in degree $d$, i.e. $\mathrm{dim}_k(M_d)$.
Moreover,
if $\ell\in S_1$ is a linear form and $R=S/I$, we denote by $R_\ell$ the ring $S/[I+(\ell)]$. 

\begin{teo}[Green's Hyperplane Restriction Theorem]\label{greenthm}
Let $I$ be a homogeneous ideal in $S$, $d\in \mathbb{N}$, then
\[
	H((S/I)_\ell,d)\leq H(S/I,d)_{\langle d \rangle},
\]
where $\ell$ is a generic linear form. Moreover, equality holds when $I_d$ is a $lex$-segment space.
\end{teo}

\section{Extension to modules}
In  the case of the restriction of an $S$-module $M$ to a generic hyperplane $\ell$
we expect  not only an upper bound on the Hilbert function $H(M_\ell,d)$, but also that this upper bound is fully 
described by lexicographic modules, in complete analogy with  Green's theorem. \\
In this section, we will be able to prove some numerical properties of the Macaulay's representations, that will lead 
to a version of  Green's HRT  for modules valid in any characteristic of the infinite field $k$.\\

\noindent
Let $F$ be a finitely generated graded free $S$-module, let us fix a homogeneous basis $\{ e_1, e_2, \dots, e_r\}$ and
let $\textrm{deg}(e_i)=f_i$,
where, without loss of generality,
we may assume that  and $f_1 \leq f_2\leq \cdots \leq f_r$.\\
We now define the monomial modules and we induce a monomial ordering on $F$ in such a way that the concept of lexicographic
module may be defined. We also define another monomial ordering on $F$, a particular $revlex$, that will let us use the tool 
of generic initial modules in our circumstances.
\begin{defn}
	A monomial in $F$ is an element of the form $me_i$ where $m\in \mathrm{Mon}(S)$. 
	A submodule $M\subseteq F$ is \textsl{monomial} if it is generated by monomials, in this 
	case it can be written as $I_1e_1 \oplus I_2e_2 \oplus \cdots \oplus I_r e_r$, where $I_i$ is a monomial ideal for $i=1,2,\dots,r$.
\end{defn}
One can extend a monomial ordering, defined in the polynomial ring $S$, to a  finitely generated free $S$-module in a really natural way,
in particular here we define the $deglex$ ordering in $F$.
\begin{defn}
	Given two monomials in $F$, $me_i$ and $ne_j$, we say that $me_i >_{deglex} ne_j$ if either $i=j$ and 
	$m>_{deglex}n$ in $S$ or $i<j$. 
	In particular, we have that $e_1>e_2>\cdots >e_r$.
\end{defn}
Again, as in the polynomial ring, when we restrict to a precise homogeneous component, the two induced monomial orderings, $deglex$ and $lex$, coincide.
\begin{defn}
	A monomial  graded submodule $L$ is a \textsl{lexicographic module} if, for every degree $d$, $L_d$ is spanned by 
	the largest, with respect 
	to the lexicographic 	order (or $deglex$), $H(L,d)$ monomials.
\end{defn}
Let us now recall the definition of reverse lexicographic order on $F$ as given in (\cite[pg. 339]{eisen}), some easy results about it.
\begin{defn}
 The \emph{reverse lexicographic order} on $F$ is defined by choosing an order on the basis of $F$, say $e_1> \cdots >e_r$ and
 by setting $me_i>_{revlex}ne_j$ iff either $\mathrm{deg}(me_i)>\mathrm{deg}(ne_j)$ or the degrees are the same and $m>_{revlex}n$ in $S$ 
 or $m=n$ and $i<j$.  
\end{defn}
\begin{rem}
 Notice that the $deglex$ is a POT (Position over Term) monomial order, on the contrary the $revlex$  is TOP (Term over Position)
 monomial order.
\end{rem}
\begin{defn}
 The initial module of a submodule  $M$, denoted by $\mathrm{in}(M)$, is the submodule of $F$ generated by the set $\{\mathrm{in}(m)| 
 \ m\in M\}$, 
 i.e. by all the leading terms (according to the lexicographic order or $deglex$ or $revlex$) of elements in $M$. 
\end{defn}
\begin{pro}\cite[Proposition 15.12]{eisen}\label{quoziente}
 Suppose that $F$ is a free $S$-module with basis $\{e_1, \dots , e_r\}$ and reverse lexicographic order. Let $M$ be a graded submodule.\\
 Then $\mathrm{in}(M+x_nF)=\mathrm{in}(M)+x_nF$ and $(\mathrm{in}(M):_F x_n)=\mathrm{in}(M:_F x_n)$.
\end{pro}

Let us define the concept  of \emph{generic initial module} as done in (\cite[Chapters 1 and 2]{pardue}).\\
Let $\mathrm{GL}(n)$ be the general linear group. It acts on $S$ by $(a_{ij})x_j=\sum_{i=1}^n a_{ij}x_i$. Let $\mathrm{GL}(F)$ be the group 
of the $S$-module
automorphisms of $F$. An element $\phi$ in $\mathrm{GL}(F)$ is a homogeneous automorphism and can be represented by a matrix $(t_{ij})$, 
with $t_{ij}\in S_{f_i-f_j}$,
where $\phi(e_i)=\sum_{j=1}^r t_{ij}e_j$. We also have an action of $\mathrm{GL}(n)$ on $\mathrm{GL}(F)$ given by $a\cdot \phi= a\phi a^{-1}$, 
for all $a\in \mathrm{GL}(n)$ and $\phi \in \mathrm{GL}(F)$. The two groups together act on $F$  by mean of their semi-direct product 
$G=\mathrm{GL}(n)\rtimes \mathrm{GL}(F)$. \\ Let us consider $B(n)$ the subgroup of $\mathrm{GL}(n)$ of upper triangular invertible matrices: 
from the action given above, a matrix in $B(n)$ will send $x_1$ to a multiple of itself, and $x_n$ to the linear
combination of all variables
$x_i$'s with coefficients the last column of tha matrix.
Let $B(F)$ be the subgroup of all automorphisms in $\mathrm{GL}(F)$ represented by lower 
triangular matrices: they send each $e_l$ to an $S$-linear
combination of $e_1, \dots, e_l$. Let $B=B(n)\rtimes B(F)$.\\
Pardue, in \cite{pardue}, proved the generalization to modules of the Galligo's theorem, in a more general setting than ours,
and in particular he proved the following result:
\begin{pro}
 Let $M\subset F$ be a graded module, and let $<$ be a monomial order on $F$, then there exists a Zariski open set $U\subseteq G$
 such that $\mathrm{in}_<(\phi M)$ is constant for every $\phi \in U$.  Moreover,  $\mathrm{in_<(\phi M)}$ is fixed by
  the action of  group $B$.
\end{pro}
\begin{defn}
 The monomial submodule $\mathrm{in}_<(\phi M)$, $\phi \in U$,
  is called the \emph{generic initial module} of $M$, and denoted by $\mathrm{gin_<(M)}$.
\end{defn}

The generic initial module has really nice properties in case we use, as monomial order on $F$, the $revlex$ order defined above.\\
\begin{pro}
For a generic linear form $\ell$ and a graded submodule $M$ of $F$, we have that
\begin{equation*}
 \mathrm{gin(M_\ell)}=\mathrm{gin(M)}_{x_n}.
\end{equation*}
\end{pro}
The $revlex$ $\mathrm{gin(M)}$ has the properties in $\ref{quoziente}$, and moreover, since it is B-fixed, we can work in 
generic coordinates, and choose the coordinates in such a way that $\ell=x_n$.

If $\ell \in S_1$ is a linear form, and $M$ a graded submodule of $F$, let us denote by $(F/M)_\ell$
the restriction of the module $F/M$ to $\ell$, which is equal to $F/(M+\ell F)$. Notice that $F_\ell$ is a module over $S_\ell$ and that $S_\ell$
can be thought as polynomial ring in $n-1$ variables.\\

Before the statement of the main theorem of the section, we are going to prove two propositions, regarding the behavior of 
the function $*_{ \langle d \rangle}:\mathbb{N}\rightarrow \mathbb{N}$: in particular we need to provide an upper bound to  the
function 
$a_{\langle d_1 \rangle} +b_{\langle d_2 \rangle}$, where $d_1 \geq d_2$, generalizing   \cite[Lemma 4.4]{gasha}; and then, we will 
extend the latter inequality to
the sum of $r=\mathrm{rank}(F)$ integers.\\
These two proposition will be just the numerical translations of the extension of Green's theorem, respectively in the 
$\mathrm{rank}(F)=2$ case, and in the $\mathrm{rank}(F)=r$ case.

\begin{lem}\cite[Lemma 4.4, Lemma 4.5]{gasha}\label{gasha}
	For any $a,b,d \in \mathbb{N}$, the two following inequalities hold:	$a_{\langle d \rangle}+b_{\langle d \rangle}\leq 
	(a+b)_{\langle d \rangle}$; 
	$a_{\langle d+1 \rangle}\leq a_{\langle d \rangle}.$ 
\end{lem}

\begin{rem}\label{herz}
If $a,d \in \mathbb{N}$ and $a=\binom{a_d}{d}+\cdots + \binom{a_\delta}{\delta}$ ($a_\delta \geq \delta$)
is the $d$-th Macaulay representation of $a$, then
$$(a-1)_{\langle d \rangle}=a_{\langle d \rangle}\ \textrm{ if and only if } \ a_\delta=\delta. $$
\end{rem}

\begin{pro}\label{ranktwo}
	Given $a,b \in \mathbb{N}$, $a \leq N_1=\binom{n+d_1-1}{d_1}$, $b\leq N_2=\binom{n+d_2-1}{d_2}$ and 
	$d_1\geq d_2, \ d_1,d_2 \in \mathbb{N}$, then
\begin{equation*}
	a_{\langle d_1 \rangle}+b_{\langle d_2 \rangle}\leq 
	\begin{cases}
		(a+b)_{\langle d_2 \rangle} &  \ \textrm{if } \ a+b\leq N_2,\\
		(a+b-N_2)_{\langle d_1 \rangle}+(N_2)_{\langle d_2 \rangle} & \ \textrm{if} \ a+b \geq N_2.
	\end{cases}
\end{equation*}
	
\end{pro}
\begin{dimo}
We will now prove the claim by induction on $i=a+b$.\\
The case $i=0$ is trivial.
Thus, let the claim be true for all $a',b'$ such that $0<a'+b'< i$.\\
Let  $a=\binom{a_{d_1}}{d_1}+\cdots +\binom{a_{d_2}}{d_2}+\cdots + \binom{a_1}{1}$, $b=\binom{b_{d_2}}{d_2}+\cdots + \binom{b_1}{1}$, 
where some of the last top coefficients in the two Macaulay representations may be zero, and let
us denote these Macaulay representations 
 as vectors $a=(a_{d_1}, \dots, a_1)$ 
and $b=(b_{d_2}, \dots , b_1)$. \\
Let us now distinguish between two cases.
\begin{description}
 \item[First Case] If the  $d_1$-Macaulay representation of $a$ ends after  degree $d_2$, i.e. $\textrm{min}\{i |\ a_i\geq i
 \}\leq d_2$, we write $a=(a_{d_1}, \dots,a_{d_2}, \dots, a_1 )$ and $b=(b_{d_2}, \dots, b_1)$, where in both vectors the last
 entries, after degree $d_2$, may be zero.\\ Then, we redefine $a$ and $b$ as follows, sorting them entry by entry:
\begin{eqnarray*}
&{}&a=(a_{d_1},\dots, a_{d_2+1}, min\{a_{d_2},b_{d_2}\},\dots ,min\{a_1,b_1\}),\\
&{}& b=(max\{a_{d_2},b_{d_2}\}, \dots ,max\{a_1,b_1\}).
\end{eqnarray*}
Observe that, in this process, we are decreasing $a$ of the same amount we increase $b$ (we are not modifying $a+b$) 
and we are not changing $a_{\langle d_1 \rangle}+b_{\langle d_2 \rangle}$, since the sum of the two binomial expansions will 
still contain the same terms. Moreover, in this way, we obtain that the Macaulay 
representation of $a$ ends sooner than that of $b$.\\
The idea is to find $a',b'$ such that $a'+b'=a+b$ and $a\leq a'$, $b' \geq b$ and such that the following are satisfied:
\begin{enumerate}
 \item $a'_{\langle d_1 \rangle}+b'_{\langle d_2 \rangle}\geq a_{\langle d_1 \rangle}+b_{\langle d_2 \rangle}$;
 \item if $i<N_2$, either $a'=0$ and $b'=i$; if $i\geq N_2$, either $a'=i-N_2$ and $b'=N_2$; 
 or at least one between $a'$ and $b'$ has Macaulay representation ending with a 
 binomial coefficient of the form $\binom{t}{t}$ (since, in this case, we can use Remark \ref{herz} and apply induction hypothesis).
\end{enumerate}
So, we can suppose that 
\begin{eqnarray*}
 a&=&(a_{d_1},\dots,a_{d_2}, a_{d_2-1}, \dots , a_s, 0, \dots, 0)\\
 b&=&(b_{d_2}, b_{d_2-1}, \dots , b_t, 0, \dots, 0),
\end{eqnarray*} 
where $a_s>s$ and $b_t>t$ (otherwise we can use Remark \ref{herz} and apply induction hypothesis), and $s\geq t$.\\
We need to show how to find $a'$ and $b'$, by decreasing $a$ and increasing $b$ of the same amount
in  three different cases.\\
In the case in which $s=t=1$, we have $a_1>1$ and  $b_1>1$,  we can consider
$(a_{d_1},\dots, a_{d_2}, \dots, a_1-1)$ and $(b_{d_2}, \dots , b_2, b_1+1)$, without
changing $a_{\langle d_1 \rangle}+b_{\langle d_2 \rangle}$. 
Applying this several times we will get to the point in which $a_1=1$, so we have reached the desired $a'$ and $b'$.\\
The second case is when $t=1$ and $s>1$, here we have to write
$\binom{a_s}{s}=\sum_{j=0}^s\binom{a_s-(j+1)}{s-j}$, and then consider
$(a_{d_1},\dots,a_{d_2}, \dots, a_{s+1}, a_s-1, \dots, a_s-s)$ and \- $(b_d, \dots, b_2,b_1+1)$ 
(we are decreasing $a$ and increasing $b$ by $1$). After this step, we are again in the first case.\\
Notice that in these first two cases if $b_1+1=b_2$ we have to
consider $(b_{d_2}, \dots, b_2+1,0)$ and then we sort again by taking the $min\{ a_i, b_i\}$ and the
$max\{ a_i, b_i\}$.\\
The last case is when $s>t>1$: here we write $\binom{a_s}{s}=\sum_{j=0}^{s-t}\binom{a_s-(j+1)}{s-j}+\binom{a_s-(s-t+1)}{t-1}$;
so we consider:
\begin{eqnarray*}
 &{}&(a_{d_1},\dots,a_{d_2}, \dots , a_{s+1}, a_s-1, \dots, a_s-(s-t+1)\ (\textrm{position $t$}),\\
 &{}& a_s-(s-t+1)\ (\textrm{position $t-1$}), 0,\dots,0)\\
 &{}&(b_{d_2}, b_{d_2-1}, \dots , b_t, 0,\dots,0),
\end{eqnarray*}
and we resort by taking the min and the max in each modified degree.\\

\item[Second Case] When the Macaulay representation of $a$ ends before the degree $d_2$, say in degree $j$, say $k$ is the 
smallest degree for which $b_k\geq k$. If $d_2>1$, then let us consider the decomposition:
$$
	\binom{a_j}{j}=\sum_{i=k+1}^j \binom{a_j-(j-i+1)}{i}+\binom{a_j-(j-k)}{k};
$$
so we write 
\begin{eqnarray*}
 a&=&(a_{d_1},\dots, a_{j+1}, a_j-1, \dots, a_j-(j-d_2+1), \dots,\\
 &{}& a_{j}-(j-k),a_{j}-(j-k), 0,\dots,0),
\end{eqnarray*}
and we compare, entry by entry, with $b$ moving the smaller entries in the new vector $a$ and the larger entries in $b$.\\
The only case in which this cannot work is in the trivial case $b=(\beta)$, $a=(a_{d_1},\dots,\alpha, \alpha)$ 
(after having used the previous decomposition) and $\alpha=\beta$, in this situation we just consider
$(a_{d_1},\dots, \alpha, \alpha-1)$ and $(\beta+1)$.\\
\end{description}
So, for all the possible cases we have found a way to subtract to $a$ and add to $b$ the same number in such a way that 
the function $f(x,y)=x_{\langle d_1 \rangle}+y_{\langle d_2 \rangle}$ does not change: we will continue, recursively, to do this until 
we get the extremal cases, or until one of the two numbers has Macaulay representation ending with $\binom{t}{t}$.

\end{dimo}

\begin{oss}
 Notice that the previous proposition is just the numerical translation of the generalization of Green's HRT 
 to the case in which we have  $F=S e_1\oplus S e_2$, with $e_1 >_{deglex} e_2$, $deg(e_1)=f_1$, $deg(e_2)=f_2$ and $f_1\leq f_2$. 
In fact, for any degree $d\in \mathbb{N}$, we set  $d_1=d-f_1$, $d_2=d-f_2$, and
$N_1=\binom{n+d_1-1}{d_1}=\textrm{dim}_k(S_{d_1})$, $N_2=\binom{n+d_2-1}{d_2}=\textrm{dim}_k(S_{d_2})$. 
\end{oss}

\begin{rem}
The previous proposition implies the first part of Lemma \ref{gasha}: in fact, when we have $f_1=f_2$, and so $N_1=N_2=N$,
we can repeat the same argument of the proof of Proposition \ref{ranktwo}, with the only difference that we continue decreasing 
$a$ and increasing $b$ of the same amount until $a=0$ and $b=i$, i.e. without distinguish between the cases in which $i\leq N$ 
and $i\geq N$.
\end{rem}
The following proposition just extends the previous inequality to the case of more summands.
\begin{pro}\label{higher2}
Let $d_1, \dots, d_r\in \mathbb{N}$ a non-increasing sequence, and set $N_i=\binom{n+d_i-1}{d_i}$. Let $a_i\leq N_i$ be an integer.
Then, if, for some $j$, $\sum_{i=j+1}^r N_i\leq \sum_{i=1}^r a_i \leq \sum_{i=j}^r N_i$, the following inequality holds:
\[
	\sum_{i=1}^r {a_i}_{\langle d_i \rangle}\leq \left(\sum_{i=1}^r a_i-\sum_{i=j+1}^r N_i\right)_{\langle d_j \rangle}+
	\sum_{i=j+1}^r {N_i}_{\langle d_i \rangle}.
\] 
\end{pro}
	\begin{dimo}
	By induction on $r$. The case $r=2$ is given by the previous result, let us now suppose that the 
	inequality holds for $s<r$.
	Let $\sum_{i=j+1}^r N_i\leq \sum_{i=1}^r a_i \leq \sum_{i=j}^r N_i$. We then distinguish between two cases.\\
	In case $a_j+a_{j+1}\leq N_{j+1}$, then ${a_j}_{\langle d_j \rangle}+{a_{j+1}}_{\langle d_{j+1}\rangle}\leq(a_j+a_{j+1})_{\langle d_{j+1} \rangle}$. \\
	Since $\sum_{i=j+1}^r N_i\leq \sum_{i=1}^r a_i\leq \sum_{i=j}^r N_i \leq N_{j-1}+\sum_{i=j+1}^r N_i$, we can apply 
	the induction to have the following:
	\begin{eqnarray*}
	 	&{}& \sum_{i=1}^r {a_i}_{\langle d_i \rangle}\leq \sum_{i=1,i\neq j,j+1}^r {a_i}_{\langle d_i \rangle}+(a_j+a_{j+1})_{\langle d_{j+1} \rangle}\leq \left(\sum_{i=1}^r a_i-\sum_{i=j+1}^r N_i\right)_{\langle d_{j-1} \rangle}+ \\
	 	 &{}& +\sum_{i=j+1}^r {N_i}_{\langle d_i \rangle}\leq  \left(\sum_{i=1}^r a_i-\sum_{i=j+1}^r N_i\right)_{\langle d_j \rangle}+\sum_{i=j+1}^r {N_i}_{\langle d_i \rangle}.
	\end{eqnarray*}
	In the case in which $a_j+a_{j+1}\geq N_{j+1}$, then, by the base of induction, ${a_j}_{\langle d_{j} \rangle}+{a_{j+1}}_{\langle d_{j+1} 
	\rangle}\leq(a_j+a_{j+1}-N_{j+1})_{\langle d_{j} \rangle}+{N_{j+1}}_{\langle d_{j+1} \rangle}$. We have that
	$a_j+a_{j+1}-N_{j+1}\leq N_j$ and $\sum_{i=j+2}^r N_i\leq \sum_{i=1}^r a_i -N_{j+1}\leq N_j+\sum_{i=j+2}^r N_i$, 
	thus by induction:
	\[
	 	\sum_{i=1,i\neq j,j+1}^r {a_i}_{\langle d_i \rangle}+(a_j+a_{j+1}-N_{j+1})_{\langle d_{j} \rangle} \leq 
	 	\left(\sum_{i=1}^r a_i-\sum_{i=j+1}^r N_i\right)_{\langle d_{j} \rangle}+\sum_{i=j+2}^r {N_i}_{\langle d_i \rangle},
	\]
	By adding ${N_{j+1}}_{\langle d_{j+1} \rangle}$ to both sides, we get 
	\[
	\sum_{i=1}^r {a_i}_{\langle d_i \rangle}\leq \left(\sum_{i=1}^r a_i-\sum_{i=j+1}^r N_i\right)_{\langle d_j \rangle}+\sum_{i=j+1}^r {N_i}_{\langle d_i \rangle}.
\] 
	\end{dimo}
	
	The idea for the extension of the  Green's theorem to modules is to reduce to the monomial case: 
	a monomial submodule is direct sum of
monomial components, 
and so we may just apply Theorem \ref{greenthm} to each of these components. 
Afterwards we need to  bound the Hilbert function of the quotient by the monomial submodule by the Hilbert function 
of the quotient by lexicographic module.\\

\noindent
Let us now give a definition that will be used inside the statement of our generalization of Green's theorem, and then we are ready
to prove our main result.
\begin{defn}
Let $M$ be a submodule in $F$, and $m\in \mathbb{N}$. Set $d_i=m-f_i $ ($\{ d_1, d_2,\dots, d_r\}$ is a non-increasing sequence) and $N_i=\binom{n+d_i-1}{d_i}= \mathrm{dim}_k S_{d_i}$.
Then, if $\sum_{i=j+1}^r N_i\leq H(F/M,m)\leq \sum_{i=j}^r N_i$, for some $j$, we define
\[
	H(F/M,m)_{\{ m, r\}}=\left(H(F/M,m)-\sum_{i=j+1}^r N_i\right)_{\langle d_j \rangle}+\sum_{i=j+1}^r {N_i}_{\langle d_i \rangle}.
\]
\end{defn}
\begin{teo}\label{main}
Let $F=Se_1 \oplus \cdots \oplus Se_r $ where $deg(e_i)=f_i$ for all $i$. Let $M$ be a submodule in $F$, then
$$
	H((F/M)_\ell,m)\leq H((F/L)_\ell,m)
$$
where $\ell$ is generic linear form, $m\in \mathbb{N}$, and $L$ is a submodule that in degree $m$ is generated by
a $lex$-segment of length $H(M,m)$.\\
Moreover, 
$$H((F/L)_\ell,m)=H(F/M,m)_{\{ m, r\}}.$$
\end{teo}
\begin{dimo}
Assume that $M$ is a monomial submodule of $F$. Otherwise, let us take the generic initial module of $M$, $\mathrm{gin}(M)$, with
respect  to  the monomial order $revlex$. We get a new submodule
which has the same Hilbert function of $M$ and moreover is a monomial submodule of $F$. 
We have that $\mathrm{gin}(M)_{x_n}=\mathrm{gin}(M_\ell)$. 
If we work in generic coordinates, we can choose $\ell=x_n$ and so we get that $H((F/M)_\ell,d)=H((F/\mathrm{gin}(M))_\ell,d)$. So we may substitute 
$M$ by $\mathrm{gin}(M)$ in order to get a monomial submodule.
Thus, we can write $M$ as $I_1e_1 \oplus \cdots \oplus I_r e_r$, where each $I_i$ is a monomial ideal, 
hence a homogeneous ideal in $R$.
\\
Let $d \in \mathbb{N}$ and $\ell$ be a generic linear form, then
\begin{eqnarray*}
	&{}&(F/M)_\ell= F/(M+\ell F) \cong F/[(I_1+(\ell))e_1 \oplus \cdots \oplus (I_r+(\ell))]\cong \\
	&{}&  (S/I_1)_\ell e_1 \oplus \cdots \oplus (S/I_r)_\ell e_r.
\end{eqnarray*}
So we have that, by applying  Green's HRT to each $(R/I_i)_\ell$,
\begin{equation}\label{ine}
	 H((F/M)_\ell,m)= \sum_{i=1}^r H((R/I_i)_\ell,m-f_i)\leq \sum_{i=1}^r H(R/I_i, m-f_i)_{\langle m-f_i \rangle}.
\end{equation}
Finally the last sum in (\ref{ine}) is bounded by $H((F/L)_\ell,m)$ where $L$ is a submodule such that $L_m$ is 
generated by a $lex$-segment of length $H(M,m)$: this follows directly from Proposition \ref{higher2}.

\end{dimo}
	
\subsection{Green's HRT after scaling for modules}
	We now derive from last theorem a result, that will give the possibility to apply Green's HRT in case it is 
	difficult to have the Macaulay binomial representations. The idea is to consider the Hilbert function up to a multiplication by
	an integer. In this way, we will obtain a linear inequality, that will hold for all modules.
	\begin{cor}
		Let $M$ be a module over $S=k[x_1, \dots, x_n]$, generated in degree $0$, and let $\ell$ be a generic linear form,
		$d\in \mathbb{N}$, then
		$$
			H(M_\ell,d)\leq \frac{n-1}{n+d-1}\ H(M,d).
		$$
	\end{cor}
		\begin{dimo}
			For any integer $\delta$, we have that $ H\left(\bigoplus \limits_{i=1}^\delta M_\ell, d\right)=\delta H(M_\ell,d)$, 
			so we can choose $\delta$ such that there is $\beta \in \mathbb{N}$ satisfying the following equality
			$$
				\delta H(M,d)=\beta \mathrm{dim}_k(S_d)=\beta \binom{n+d-1}{d},
			$$
			so by Theorem \ref{main} we have
			\begin{eqnarray*}
				 &{}& H\left(\bigoplus \limits_{i=1}^\delta M_\ell, d\right)\leq \beta  \binom{n+d-1}{d}_{\langle d \rangle}=\beta 
				 \binom{n+d-2}{d}=\\
				  &{}& \beta \binom{n+d-1}{d}\frac{n-1}{n+d-1}= \frac{n-1}{n+d-1}H\left(\bigoplus \limits_{i=1}^\delta M, d\right),
			\end{eqnarray*}
			then, dividing by $\delta$, we get the claim.
		\end{dimo}

\begin{oss}
This corollary can be described by the expression "after scaling" because the last fraction satisfies
			the following equality
			$$
			\frac{n-1}{n+d-1}=\frac{\dim(R_d)}{\dim(S_d)},
			$$
			where $S=k[x_1, x_2, \dots, x_n]$ and $R=k[x_1, x_2,\dots,x_{n-1}]\cong S_\ell$.\\
			Therefore, the claim of the Corollary may be written as:
			$$
				\frac{H(M_l,d)}{H(S_l,d)} \leq \frac{H(M,d)}{H(S,d)},
			$$
			so we have scaled by respective value of the Hilbert function of the rings over which we have the two modules.
		\end{oss}
		
\section{Application to level algebras}		
In this section we will apply the extension of the Green's HRT to level algebras.
\begin{defn}
 Let $A=S/I=A_0\oplus \cdots \oplus A_c$ be an artinian $k$-algebra. Let us consider the 
 \emph{socle} of $A$, which is defined by $\mathrm{soc}(A)=\mathrm{ann}_\mathbf{m}(A)$,
 where $\mathbf{m}=\oplus_{i>0}A_i$. If $\mathrm{soc}(A)=\mathcal{U}_0\oplus \cdots \oplus 
 \mathcal{U}_c$, we say that $A$ is \emph{level of type $b$} if $\mathrm{soc}(A)=\mathcal{U}_c=A_c$ and $\mathrm{dim}(A_c)=b$. 
 We say that $A$ is \emph{Gorenstein} when it is level of type $1$.
\end{defn}

\noindent
Let $A=A_0\oplus \cdots \oplus A_c$ be a level algebra with Hilbert function $(h_0,h_1,\dots,h_c)$.
By considering the dual of $A$, $A^\vee=\mathrm{Hom}_k(A,k)$, we get a \emph{level module}, since
is an artinian $S$-module with generators in a single degree, and with one-dimensional socle, i.e. also concentrated in a single degree.\\
We want to consider the multiplication by a generic linear form $\ell \in A_1$.\\
We first apply Theorem \ref{greenthm}. In this way we get the following bound on the Hilbert function 
of the ideal $(\ell)$.
$$
  \mathrm{dim}_k (\ell)_i \geq h_i-(h_i)_{\langle i \rangle}, \ 0\leq i \leq c.
$$
We can also apply Theorem \ref{main} to  the level module $M=A^\vee(-c)$. Notice that $M$
has Hilbert function $(h_c,h_{c-1},\dots ,h_0)$. We get in this case the following bound
for the module $\ell M$:
$$
\mathrm{dim}(\ell M)_i \geq h_{c-i}-(h_{c-i})_{\{ i\}}, \ 0\leq i \leq c.
$$
In the right hand side of this inequality,  $a_{\{ i\}}$ denotes $q(s_i)_{\langle  i\rangle}+
r_{\langle  i\rangle}$, where $a=qs_i+r$ and $s_i=\binom{i+n-1}{i}$. 

We observe that $(\ell) \cong (\ell M)^\vee(-c+1)$. This implies that  we can give two lower bounds for the Hilbert function of $(\ell)$, namely $h^{GM}_i$ and $h^{G}_i$, where we set  $h^{GM}_i=h_i-(h_i)_{\{c-i \}}$ and $h^{G}_i=h_{i+1}-(h_{i+1})_{\langle i+1 \rangle}$.

We would like to find conditions in which $h^{GM}_i> h^G_i$, so when the bound given by
the HRT for modules is better than the one given by the Green's HRT.\\

In \cite{Level}, Geramita, Harima, Migliore and Shin, in case the polynomial ring $S$ has three variables, provided a list of all the possible Hilbert functions for level algebras of socle in degree $3,4,5$ and in degree $6$ and type $2$.
Among these, the condition $h^{GM}_i> h^G_i$ is true only  for the ones contained in Table \ref{table:comparison} in page \pageref{table:comparison}.

\begin{pro}
 For $i$ such that the following conditions are all verified, we have that $h^{GM}_i\geq  h^G_i$:
    \begin{eqnarray*}
    \begin{cases}
        i-1\leq c-i;\\
	h_i=h_{i+1};\\
	s_{c-i}> h_{i+1}.
	\end{cases}
    \end{eqnarray*}
\end{pro}
\begin{dimo}
 We have that $h^{GM}_i=h_i-(h_i)_{\{ c-i \}}=h_{i+1}-(h_{i+1})_{\{ c-i \}}$.
 Since $s_{c-i}> h_{i+1}$, $(h_{i+1})_{\{ c-i \}}=(h_{i+1})_{\langle c-i \rangle} $.
 We may then compare $(h_{i+1})_{\langle c-i \rangle}$ and $(h_{i+1})_{\langle i+1 \rangle}$, using 
 Lemma \ref{gasha}, and the fact that $i-1\leq c-i$, so we get that $h^{GM}_i\geq  h^G_i$.
\end{dimo}

The previous condition does not characterize all cases in which $h^{GM}_i> h^G_i$, but only the ones in bold font in the table.  Notice that  if in the first column we have a certain number $i$ this means, according to our notation, that $h^{GM}_{i-1}>h^{G}_{i-1}$.\\

\begin{que}
What are the sufficient conditions for which the bound given by the HRT for modules is better than the bound given by the classical HRT?
\end{que}

\begin{table}[h]
\begin{tabular}{| l | l | l | l |}

\hline
&&& \\
 Position & $h$ & $h^{GM}$ & $h^{G}$\\
\hline
 2 & \textbf{1, 3, 3, 3, 2} & 1, 3, 2, 1 & 1, 2, 3, 2\\
\hline
2 & \textbf{1, 3, 3, 3, 3 }& 1, 3, 2, 1 & 1, 2, 3, 3\\
\hline
2 & \textbf{1, 3, 3, 3, 3, 2} & 1, 3, 3, 2, 1 & 1, 2, 3, 3, 2\\
\hline
4 & 1, 3, 6, 8, 5, 2 & 1, 3, 5, 5, 2 & 1, 3, 6, 4, 2\\
\hline
4 & 1, 3, 6, 9, 5, 2 & 1, 3, 5, 5, 2 & 1, 3, 6, 4, 2 \\
\hline
4 & 1, 3, 6, 10, 6, 2 & 1, 3, 5, 6, 2 & 1, 3, 6, 5, 2\\
\hline
2 & \textbf{1, 3, 3, 3, 3, 3} & 1, 3, 3, 2, 1 & 1, 2, 3, 3, 3\\
\hline
4 & 1, 3, 6, 8, 5, 3 & 1, 3, 5, 5, 2 & 1, 3, 6, 4, 3\\
\hline
4 & 1, 3, 6, 10, 6, 3 & 1, 3, 5, 6, 2 & 1, 3, 6, 5, 3\\
\hline
2 & \textbf{1, 3, 3, 3, 3, 3, 2} & 1, 3, 3, 3, 2, 1 & 1, 2, 3, 3, 3, 2\\
\hline
3 & \textbf{1, 3, 4, 4, 4, 3, 2} & 1, 3, 4, 3, 3, 1 & 1, 3, 3, 4, 3, 2\\
\hline
3 & \textbf{1, 3, 4, 4, 4, 4, 2} & 1, 3, 4, 3, 3, 2 & 1, 3, 3, 4, 4, 2\\
\hline
4 & 1, 3, 5, 6, 5, 4, 2 & 1, 3, 4, 5, 3, 2 & 1, 3, 5, 4, 4, 2\\
\hline
4& 1, 3, 6, 8, 6, 4, 2 & 1, 3, 5, 6, 3, 2 & 1, 3, 6, 5, 4, 2\\
\hline
5 & 1, 3, 6, 8, 10, 6, 2 & 1, 3, 5, 6, 6, 2 & 1, 3, 6, 8, 5, 2\\
\hline
5 & 1, 3, 6, 9, 10, 6, 2 & 1, 3, 5, 6, 6, 2 & 1, 3, 6, 8, 5, 2\\
\hline
5 & 1, 3, 6, 9, 11, 6, 2 & 1, 3, 5, 6, 6, 2 & 1, 3, 6, 9, 5, 2\\
\hline
5 & 1, 3, 6, 9, 12, 6, 2 & 1, 3, 5, 6, 6, 2 & 1, 3, 6, 9, 5, 2\\
\hline
5 & 1, 3, 6, 10, 10, 6, 2 & 1, 3, 5, 6, 6, 2 & 1, 3, 6, 8, 5, 2\\
\hline
5& 1, 3, 6, 10, 11, 6, 2 & 1, 3, 5, 6, 6, 2 & 1, 3, 6, 9, 5, 2\\
\hline
5 & 1, 3, 6, 10, 12, 6, 2 & 1, 3, 5, 6, 6, 2 & 1, 3, 6, 9, 5, 2\\
\hline 
\end{tabular}
\caption{Hilbert functions of Level artinian algebras}
\label{table:comparison}
\end{table}

\noindent
{\bf Acknowledgments.} I would like to thank my supervisor Mats Boij  for his great support, and also  Ralf Fr\"oberg for the useful comments. Moreover, I also thank Giulio Caviglia for his suggestions, that were crucial in the last phase of preparation of the paper.

\addcontentsline{toc}{section}{Bibliography}
\bibliographystyle{siam}
\bibliography{Green}

\noindent
 {\scshape Royal Institute of Technology, Department of Mathematics, S-10044 Stockholm, Sweden}.\\
 {\itshape E-mail address}: \texttt{ogreco@kth.se}

\end{document}